\documentclass{gtpart}

\agtart

\usepackage{amsxtra}
\usepackage{amsfonts}
\usepackage{amscd}
\usepackage{newlfont}
\usepackage[dvips]{graphicx}
\usepackage[all]{xy}

\title{Rational blow-down along Wahl type plumbing trees of spheres}

\author{Maria Michalogiorgaki}
\givenname{Maria}
\surname{Michalogiorgaki}
\address{Department of Mathematics\\
Princeton University\\\newline
Fine Hall\\Washington Road\\Princeton
NJ 08544\\USA}
\email{mariam@math.princeton.edu}



\volumenumber{}
\issuenumber{}
\publicationyear{}
\papernumber{}
\startpage{}
\endpage{}
\doi{}
\MR{}
\Zbl{}
\received{}
\revised{}
\accepted{}
\published{}
\publishedonline{}
\proposed{}
\seconded{}
\corresponding{}
\editor{}
\version{}

\newtheorem{theorem}{Theorem}
\newtheorem{prop}{Proposition}
\newtheorem{corollary}{Corollary}
\newtheorem{remark}{Remark}
\newtheorem{note}{Note}
\newtheorem{definition}{Definition}

\begin{document}

\begin{abstract}
In this article, we construct smooth 4--manifolds homeomorphic but
not diffeomorphic to $ \mathbb{CP}^{2} \sharp k
\overline{\mathbb{CP} ^{2}} $, for $k \in \{6,7,8,9\},$ using the
technique of rational blow-down along Wahl type plumbing trees of
spheres. (see \cite{Wa})
\end{abstract}

\maketitle

\def\co{\colon\thinspace}

\section{Introduction}

Over the past three years, and due to examples constructed by
J. Park, R. Fintushel, R. Stern, A. Stipsicz and Z. Szab\'{o} (see
\cite{Pa}, \cite{Oz&Sz}, \cite{St&SZ2}, \cite{FiSt} and
\cite{St&Sz3}), there has been renewed interest in the problem of
finding the smallest k for which $\mathbb{CP}^{2} \sharp k
\overline{\mathbb{CP} ^{2}} $ admits an exotic smooth structure. All
these examples are constructed using the rational blow-down
operation along lens spaces.

In this paper, we study a more generalized rational blow-down
operation along certain Seifert fibered 3--manifolds. This
technique, for the case of Wahl type plumbing trees of spheres,
together with knot surgery along a regular fiber in a double node
neighborhood (see \cite{FiSt}), are then used to construct manifolds
homeomorphic but not diffeomorphic to $ \mathbb{CP}^{2} \sharp k
\overline{\mathbb{CP} ^{2}} $, for $k \in \{6,7,8,9\}$.

\begin{flushleft}
$\mathbf{Aknowledgements}$
\end{flushleft}
The author would like to thank Zolt\'{a}n Szab\'{o} for his guidance
and support while working on this problem. The author is indebted to
him for sharing his ideas on the subject. The author would also like
to thank Andr\'{a}s Stipsicz for useful suggestions on an earlier
version of this article as well as the referee for a careful reading
of the manuscript and several constructive comments.

\section{Seiberg--Witten invariants and surgery along monopole L--Spaces}

In this section, we provide a very brief review of the
Seiberg--Witten theory of 4--manifolds in general as well as in the
special case $b_{2}^{+}=1$ and we study the effect of surgery along
monopole L--spaces on Seiberg--Witten invariants, using monopole
Floer homology. These will be the main tools for our constructions
in the next sections. For more details, we refer the reader to
\cite{Mo}, \cite{Sc}, \cite{FS}, \cite{KMOS} and
\cite{KM}.\\

\subsection{Seiberg--Witten invariants}

Let $X$ be an oriented, closed, Riemannian 4--manifold and
$\mathfrak{s}$ a $spin^{c}$ structure on $X$. Suppose that
$W_{\mathfrak{s}}^{+}$ and $W_{\mathfrak{s}}^{-}$ are the associated
U(2) spinor bundles and $L \to X$ with $L \simeq
detW_{\mathfrak{s}}^{+} \simeq detW_{\mathfrak{s}}^{-}$ is the
associated determinant line bundle. Given a pair $(A,\Psi) \in
A_{X}(L)\times \Gamma(W_{\mathfrak{s}}^{+})$, where $A_{X}(L)$
denotes the space of connections on L, and a $g$-self-dual 2-form $
\eta \in \Omega_{g}^{+}(X,\mathbb{R})$, the \textbf{\emph{perturbed
Seiberg--Witten equations}} are
\begin{equation}\label{eq:SW}
D_{A} \Psi = 0, \qquad F_{A}^{+}=i(\Psi\otimes\Psi^{*})_{o}+i \eta
\end{equation}
where $D_{A}\co \Gamma(W_{\mathfrak{s}}^{+}) \to
\Gamma(W_{\mathfrak{s}}^{-})$ is the Dirac operator and $ (\Psi
\otimes \Psi^{*})_{o}$ is the trace free part
of the ednomorphism $\Psi \otimes \Psi^{*}$. \\

The quotient of the solution space to the equations above under the
action of the gauge group Aut(L)=Map(X,$S^{1}$), denoted here by
$M_{X}(L)$, has \textbf{\emph{formal dimension}}
\begin{equation}\label{eq:dim}
dim M_{X}(L)= \frac{1}{4}(c_{1}(L)^{2}-(3sign(X)+2e(X)))
\end{equation}
Under the additional assumption that $b_{2}^{+}>0$ and for generic
form $\eta$, $M_{X}(L)$ is a smooth compact manifold, since when
$b_{2}^{+}>0$ there are no reducible solutions, i.e. no
singularities in the quotient space. \\

The \textbf{\emph{Seiberg--Witten invariant}} for $X$ is a function
$SW_{X}\co Spin^{c}(X) \to \mathbb{Z}$ defined as follows:
\begin{itemize}
 \item If $dimM_{X}(L)<0$ or odd, then $SW_{X}(L)=0.$
 \item If $dimM_{X}(L)=0$, then $SW_{X}(L)= $number of points in
$M_{X}(L)$, counted with sign.
 \item  If $dimM_{X}(L)=2n>0$, then
$SW_{X}(L)$$=<\mu^{n},[M_{X}(L)]>$, where $\mu \in H^{2}(M_{X}(L);
\mathbb{Z})$ is the Euler class of the basepoint map
$\widetilde{M_{X}}(L) = \{ solutions (A,\Psi) \} / Aut^{o}(L) \to
M_{X}(L)$, which is an $S^{1}$ fibration if there are no reducible
solutions. Here, $Aut^{o}(L)=$ \{gauge transformations which are the
identity on the fiber of $L$ over a fixed basepoint on $X$ \}.
\end{itemize}

$SW_{X}$ is independent of $g$ and $\eta$ provided that
$b_{2}^{+}>1$. In the case $b_{2}^{+}(X)=1$, there is a
codimension-one submanifold of metrics for which there are reducible
solutions and this must be excluded. Then the SW invariant of a
given $spin_{c}$ structure has two values, depending on the metric,
and the wall-crossing formula describes the relation between these
values.

\textit{\textbf{Wall-crossing formula}}: Suppose that $X$ is a
closed, oriented 4--manifold with $b_{2}^{+}(X)=1$,
$H_{1}(X;\mathbb{Z})=0$ and a fixed orientation for
$H_{+}^{2}(X;\mathbb{R})$, $\mathfrak{s}$ is a $spin^{c}$ structure
on $X$ such that $c_{1}(L)\neq 0$, $R$ is the space of Riemannian
metrics $g$ on $X$, $\omega^{+}(g)$ is the $g$-self-dual harmonic
form of norm one which lies in the positive component of
$H_{2}^{+}(X;\mathbb{R})$ as measured by the given orientation and
$R^{+}= \{g \in R / \omega^{+}(g) \cdot c_{1}(L)
> 0 \}$, $R^{-}= \{g \in R / \omega^{+}(g) \cdot c_{1}(L) < 0 \}$.
Then, $\forall g \in R^{+}\bigsqcup R^{-}$, $SW_{g}(\mathfrak{s})$
is defined and assuming that $d(\mathfrak{s})=dimM_{X}(L) \geq 0$
and even,
\begin{equation}\label{eq:Wall crossinq}
SW_{+}(\mathfrak{s})=SW_{-}(\mathfrak{s})-(-1)^{
\frac{d(\mathfrak{s})}{2}}.
\end{equation}
Here, $SW_{+(-)}(\mathfrak{s})$ denotes the constant value of
$SW_{g}(\mathfrak{s})$ on $R_{+(-)}$ respectively.

\subsection{Monopole Floer homology}
We now provide a very brief review of monopole Floer homology
as constructed by P. Kronheimer and T. Mrowka. We refer the reader
to \cite{KMOS} and \cite{KM} for more details and point out that
this version of Floer homology is conjectured to be isomorphic to
Heegaard Floer homology.

Let Y be a smooth, oriented, compact, connected 3--manifold without
boundary. To it, there are associated three vector spaces over a
field $\mathbb{F}$, namely $\check{HM}\textbf{.}(Y)$,
$\widehat{HM}\textbf{.}(Y)$ and $\overline{HM}\textbf{.}(Y)$. These
spaces, called Floer homology groups, come equipped with linear maps
$i_{*}, j_{*}$ and $p_{*}$ which form a long exact sequence
\begin{equation}\label{eq:exact sequence}
... \xrightarrow{\text{$i_{*}$}}
\check{HM}\textbf{.}(Y)\xrightarrow{\text{$j_{*}$}}
\widehat{HM}\textbf{.}(Y)\xrightarrow{\text{$p_{*}$}}
\overline{HM}\textbf{.}(Y) \xrightarrow{\text{$i_{*}$}}
\check{HM}\textbf{.}(Y) \xrightarrow{\text{$j_{*}$}} ...
\end{equation}
and with an endomorphism $u$ of degree -2 that makes the three
spaces modules over the polynomial ring $\mathbb{F}[u]$.

In addition, to each cobordism $W\co Y_{0} \to Y_{1}$, there are
associated maps $\check{HM}(W)\co \check{HM}\textbf{.}(Y_{0}) \to
\check{HM}\textbf{.}(Y_{1})$, $\widehat{HM}(W)\co
\widehat{HM}\textbf{.}(Y_{0}) \to \widehat{HM}\textbf{.}(Y_{1})$ and
$\overline{HM}(W)\co \overline{HM}\textbf{.}(Y_{0}) \to
\overline{HM}\textbf{.}(Y_{1})$ for which $i_{*},j_{*},p_{*}$ give
natural transformations. These maps respect the module structure of
the Floer groups.

In this setting, we have the following
\begin{definition}
\label{def:L--space}
A rational homology 3--sphere Y for which
$j_{*}\co \check{HM}\textbf{.}(Y) \rightarrow
\widehat{HM}\textbf{.}(Y)$ is trivial is called a \textbf{monopole
L--space}.
\end{definition}

\subsection{Using monopole Floer homology to compute Seiberg--Witten invariants after
surgery along monopole L--spaces.}

After recalling some basics of Seiberg--Witten theory and
Monopole Floer homology, we move on to compute how SW invariants
change under surgery along monopole L--spaces.

Suppose X is a 4--manifold decomposed into two pieces Z and P along
a monopole L--space Y, with P negative definite, and $ \mathfrak{s}
\in Spin^{c}(X)$ (See Figure \ref{1}). Consider B another negative
definite 4--manifold bounded by Y such that $\mathfrak{s}|_{Z}$
extends to B and replace P with B in X to get $X'=Z\bigcup_{Y}B$
(See Figure \ref{2}). Denote the $spin^{c}$ structure on $X'$ by
$\mathfrak{s}'$. We would like to compute the change in SW
invariants after such an operation. To this end, we study these
configurations, using properties of monopole Floer homology. Most of
the statements we make here are discussed in \cite{KM} and
\cite{KMOS}.

\begin{remark}
\label{rem:1}
Note that our constructions in this paper (see next
section) are special cases of the above, for B rational balls and P
Wahl-type plumbing trees of spheres. The boundaries of the latter
were proven to be monopole L--spaces in \cite{KMOS}.
\end{remark}

\begin{figure}[h]
\centering
\input{manifolds.pstex_t}
\caption{$X=Z\bigcup_{Y}P$} \label{1}
\end{figure}

\begin{figure}[h]
\centering
\input{manifolds2.pstex_t}
\caption{$X'=Z\bigcup_{Y}B$} \label{2}
\end{figure}

Our goal in this subsection is to prove that

\pagebreak

\begin{theorem}
\label{thm:main}
Suppose $Y$ is a monopole L--space, $H_{1}(Y)$ is
finite, $P,B$ are negative definite 4--manifolds with $b_{1}=0$ and
$X= Z \bigcup_{Y} P$, $X' = Z \bigcup_{Y} B $, for some 4--manifold
Z. If $\mathfrak{s} \in Spin^{c}(X)$, $\mathfrak{s}' \in
Spin^{c}(X')$, $d(\mathfrak{s}), d(\mathfrak{s}') \geq 0$ and
$\mathfrak{s} \vert
_{Z} = \mathfrak{s}' \vert _{Z}$, then $SW_{X}(\mathfrak{s})=SW_{X'}(\mathfrak{s}').$ \\
In the case $b_{2}^{+}(X)=1$,
$SW_{X,a_{1}}(\mathfrak{s})=SW_{X',a_{2}}(\mathfrak{s}')$, where
$a_{1} \in H_{2}(X; \mathbb{Z}), a_{2} \in H_{2}(X'; \mathbb{Z})$,
$a_{1} \vert _{P} = a_{2} \vert _{B} = 0 $ and $ a_{1} \vert _{Z}
= a_{2} \vert_{Z}$.
\end{theorem}

\begin{proof}
Denote $W, W_{1},W_{2}$ the cobordisms $Z-B^{4}\co
S^{3} \to Y$, $P-B^{4}\co Y \to S^{3}$, $B-B^{4}\co Y \to S^{3}$
respectively and $\mathfrak{s}_{1}=\mathfrak{s}\vert _{W_{1}}$,
$\mathfrak{s}_{2}=\mathfrak{s'}\vert_{W_{2}}$.

According to Proposition 2.6 of \cite{KMOS}, the fact that Y is a
rational homology sphere implies that
\begin{equation}\label{HM bar}
\overline{HM}\textbf{.}(Y,\mathfrak{s}\vert_{Y})\simeq
\mathbb{F}[u^{-1},u]] \ as \ topological \ F[[u]] \ modules.
\end{equation}
Here $\mathbb{F}[u^{-1},u]]$ denotes Laurent series finite in the
negative direction. In addition, the long exact sequence
(\ref{eq:exact sequence}) gives the exact sequence
\begin{equation}\label{eq:short ex. seq.}
0 \to
\widehat{HM}\textbf{.}(Y)\xrightarrow{p_{*}}\overline{HM}\textbf{.}(Y)
\xrightarrow{i_{*}}\check{HM}\textbf{.}(Y) \to 0
\end{equation}
for Y monopole L--space (see Definition \ref{def:L--space}).
Combining (\ref{HM bar}) and (\ref{eq:short ex. seq.}), we get that
the sequences
\begin{displaymath} 0 \to
\widehat{HM}\textbf{.}(Y,
\mathfrak{s}\vert_{Y})\xrightarrow{p_{*}}\overline{HM}\textbf{.}(Y,
\mathfrak{s}\vert_{Y}) \xrightarrow{i_{*}}\check{HM}\textbf{.}(Y,
\mathfrak{s}\vert_{Y}) \to 0
\end{displaymath}
and
\begin{displaymath}
0 \to \mathbb{F}[[u]] \to \mathbb{F}[u^{-1},u]] \to
\mathbb{F}[u^{-1},u]] / \mathbb{F}[[u]] \to 0
\end{displaymath}
are isomorphic as sequences of topological $\mathbb{F}[[u]]$
modules. The corresponding isomorphism of short exact sequences
holds if we consider $S^{3}$ instead of $Y$, because $S^{3}$ is a
monopole L--space as a 3--manifold with positive scalar curvature.
Such 3--manifolds were proven to be monopole L--spaces in
\cite{KMOS}.

We will now use $\overrightarrow{HM}\textbf{.}$ as defined in
\cite{KM}: For $Y_{0}, Y_{1}$ compact, connected, oriented
3--manifolds and $W$ isomorphism class of connected cobordisms
equipped with an homology orientation,
$\overrightarrow{HM}\textbf{.}(W)\co
\widehat{HM}\textbf{.}(Y_{0})\rightarrow
\check{HM}\textbf{.}(Y_{1})$ is a (canonical choice of) map such
that the diagram

\begin{displaymath}
\xymatrix{ \overline{HM}\textbf{.}(Y_{0}) \ar[r]^{i_{*}} \ar[d] &
\check{HM}\textbf{.}(Y_{0}) \ar[r]^{j_{*}} \ar[d] &
\widehat{HM}\textbf{.}(Y_{0})
  \ar[r]^{p_{*}} \ar[dl]_{\overrightarrow{HM}\textbf{.}(W)} \ar[d] & \overline{HM}\textbf{.}(Y_{0}) \ar[d]  \\
\overline{HM}\textbf{.}(Y_{1})  \ar[r]^{i_{*}} &
\check{HM}\textbf{.}(Y_{1}) \ar[r]^{j_{*}} &
\widehat{HM}\textbf{.}(Y_{1})
  \ar[r]^{p_{*}} & \overline{HM}\textbf{.}(Y_{1})
}
\end{displaymath}

\noindent commutes. For the special case of a cobordism where $W$ is
the complement of two disjoint balls in a closed, oriented manifold
$X$ viewed as a cobordism $W\co S^{3} \to S^{3}$, Proposition 3.6.1
of \cite{KM} states that the sum of the SW invariants of the
4--manifold $X$ is determined by the map
$\overrightarrow{HM}\textbf{.}(W)$. Even more, using local
coefficients, it can be proven that the individual SW invariants are
determined by the same map. (See Proposition 3.8.1 in \cite{KM})
Applying the composition laws that $\overrightarrow{HM}\textbf{.}$
satisfies according to \cite{KM} to our cobordisms $W, W_{1}$ and
$W_{2}$ gives that $\overrightarrow{HM}\textbf{.}(W_{i}\circ
W)=\check{HM}\textbf{.}(W_{i})\circ
\overrightarrow{HM}\textbf{.}(W)$, $i \in \{1,2\}$. In a more
refined version, for our fixed $spin^{c}$ structures $\mathfrak{s}$
and $\mathfrak{s_{i}}$ on $W$ and $W_{i}$ respectively, we have that
\begin{equation}\label{eq:composition law}
\sum_{\substack{\mathfrak{s}_{i}' \in Spin^{c}(W\bigcup W_{i}) \\
\mathfrak{s}_{i}' \vert _{W_{i}}= \mathfrak{s}_{i} \\
\mathfrak{s}_{i}' \vert _{W}= \mathfrak{s}}}
\overrightarrow{HM}\textbf{.}(W_{i}\circ W,
\mathfrak{s}_{i}')=\check{HM}\textbf{.}(W_{i},
\mathfrak{s}_{i})\circ \overrightarrow{HM}\textbf{.}(W,
\mathfrak{s})
\end{equation}

\noindent This sum contains precisely one term, since $Y$ is a
rational homology 3--sphere and so, (\ref{eq:composition law}) gives
that
\begin{equation}\label{eq:composition law 2}
\overrightarrow{HM}\textbf{.}(W_{1}\circ W,
\mathfrak{s})=\check{HM}\textbf{.}(W_{1}, \mathfrak{s}_{1})\circ
\overrightarrow{HM}\textbf{.}(W, \mathfrak{s})
\end{equation}

\noindent and

\begin{equation}\label{eq:composition law 3}
\overrightarrow{HM}\textbf{.}(W_{2}\circ W,
\mathfrak{s}')=\check{HM}\textbf{.}(W_{2}, \mathfrak{s}_{2})\circ
\overrightarrow{HM}\textbf{.}(W, \mathfrak{s})
\end{equation}

\noindent Taking into account that
$\overrightarrow{HM}\textbf{.}(W_{1}\circ W, \mathfrak{s})$,
$\overrightarrow{HM}\textbf{.}(W_{2}\circ W, \mathfrak{s}')$
determine $SW_{X}(\mathfrak{s})$ and $SW_{X'}(\mathfrak{s}')$
respectively, it suffices to show that the maps
$\check{HM}\textbf{.}(W_{1},\mathfrak{s}_{1})$ and
$\check{HM}\textbf{.}(W_{2},\mathfrak{s}_{2})$ are isomorphisms in
the range of $\overrightarrow{HM}\textbf{.}(W, \mathfrak{s})$ to
finish the proof.

To see this, consider the diagram

\begin{equation}\label{eq:diagram}
\begin{CD}
  ... @>>> \overline{HM}\textbf{.}(Y,\mathfrak{s}\vert_{Y})  @>i_{*}>> \check{HM}\textbf{.}(Y,\mathfrak{s}\vert_{Y}) @>>>
  0 \\
      @. @V\overline{HM}\textbf{.}(W_{i},\mathfrak{s}_{i})VV @V\check{HM}\textbf{.}(W_{i},\mathfrak{s}_{i})VV @. \\
  ... @>>> \overline{HM}\textbf{.}(S^{3},\mathfrak{s}_{i}\vert_{S^{3}})  @>i_{*}>> \check{HM}\textbf{.}(S^{3},\mathfrak{s}_{i}\vert_{S^{3}}) @>>>
  0 \\
\end{CD} \qquad i \in \{1,2\}.
\end{equation}

\noindent We have assumed that $b_{1}(W_{i})=0$ and $W_{i}$ is
negative definite and under these assumptions the map
$\overline{HM}(W_{i},\mathfrak{s}_{i})\co
\overline{HM}_{j_{0_{i}}}(Y) \to \overline{HM}_{j_{1_{i}}}(S^{3})$,
$i \in \{1,2\}$, where $j_{0_{i}} \in
J(Y,\mathfrak{s}_{i}\vert{Y})=$\{homotopy classes of oriented
2-plane fields on $Y$ that determine the $spin^{c}$ structure
$\mathfrak{s}_{i}\vert{Y}$ on $Y$\}, $j_{1_{i}} \in
J(S^{3},\mathfrak{s}_{i}\vert{S^{3}})$ and
$j_{0_{i}}\overset{\mathfrak{s}_{i}}{\sim}j_{1_{i}}$, is an
isomorphism, as was proven in \cite{KMOS}. This implies that
$\check{HM}\textbf{.}(W_{i},\mathfrak{s}_{i})$ is an isomorphism on
the range of $\overrightarrow{HM}\textbf{.}(W, \mathfrak{s})$ in the
case where $d(\mathfrak{s}),d(\mathfrak{s}')\geq 0$.

In the case $b_{2}^{+}(X)=1$, the SW invariants depend on the choice
of metric $g$ and perturbation $\eta$. If $b_{2}^{+}(W)=1$, then a
choice of $g$ and $\eta$ for $W$ determines the chamber that will be
used for our computations.
\end{proof}

\section{The topological constructions}

To construct our 4--manifolds, we will blow-down certain Wahl
type plumbing trees of spheres in rational surfaces. In order to
locate these configurations in such surfaces, we use specific
elliptic fibrations of E(1) in each case. Proofs for the existence
of such fibrations are postponed until the last section of our
article.

\subsection{An exotic smooth structure on $ \mathbb{CP}^{2} \sharp
9 \overline{\mathbb{CP} ^{2}}$.}

Our first construction relies on the following existence
result.

\begin{prop}
\label{prop:1}
There is an elliptic fibration of E(1) $ \to
\mathbb{CP}^{1}$ with a singular $I_{3}$ fiber, 9 fishtail fibers
and one section.
\end{prop}

\begin{proof} An outline of the proof of this
proposition is provided in the appendix.
\end{proof}

Consider the plumbing tree of spheres\\

\noindent \setlength{\unitlength}{3cm}
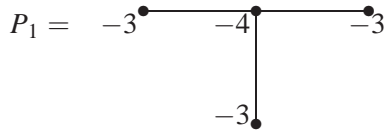
          \begin{picture}(1.5,.5)
          \put(.6,.5){\circle*{.05}}
          \put(.4,.4){$-3$}
          \put(1.1,.5){\circle*{.05}}
          \put(.9,.4){$-4$}
          \put(1.6,.5){\circle*{.05}}
          \put(1.5,.4){$-3$}
          \put(1.1,0){\circle*{.05}}
          \put(.9,0){$-3$}
          \put(.6,.5){\line(1,0){1}}
          \put(1.1,0){\line(0,1){.5}}
          \put(0,.4){$P_{1}=$}
          \end{picture}

\noindent where dots represent disk bundles over $S^{2}$, numbers
assigned to them refer to the corresponding Euler numbers and edges
stand for plumbing connections and call $Y_{1}$ the boundary of
$P_{1}$, i.e. $Y_{1}= \partial P_{1}$. Proposition \ref{prop:2}
below provides an embedding of $P_{1}$ into $\mathbb{CP}^{2} \sharp
13 \overline{\mathbb{CP} ^{2}}$.

\begin{prop}
\label{prop:2}
$P_{1}$ embeds into $\mathbb{CP}^{2} \sharp 13
\overline{\mathbb{CP} ^{2}}$.
\end{prop}

\begin{figure}[h]
\centering
\input{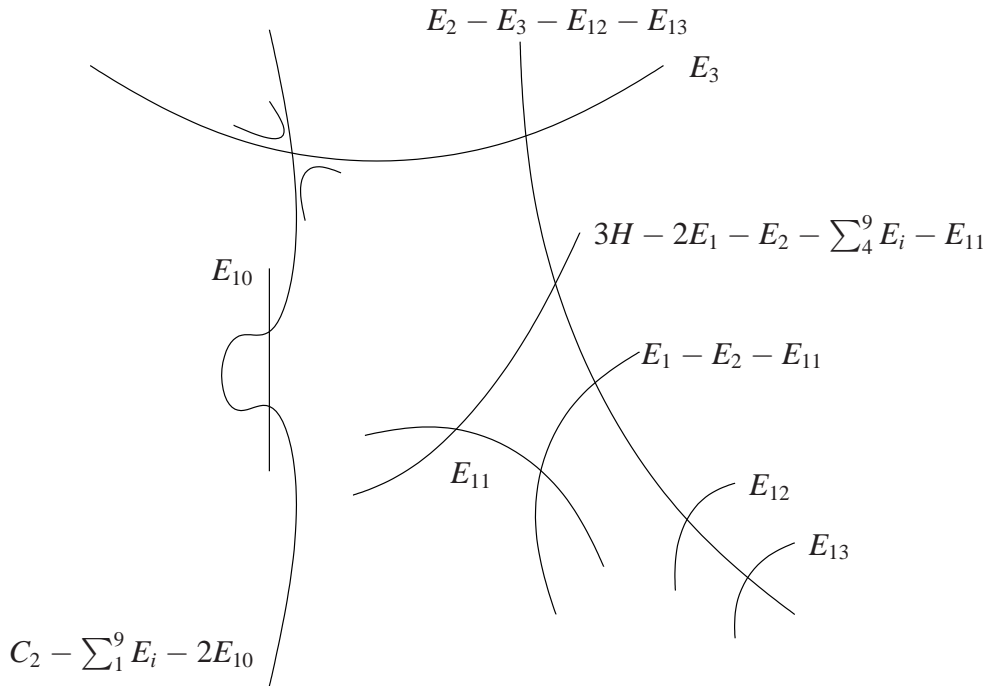}
\caption{$P_{1}$ in $\mathbb{CP}^{2} \sharp 13 \overline{\mathbb{CP}
^{2}}$} \label{3}
\end{figure}

\begin{proof} We use the fibration of Proposition \ref{prop:1} and
in particular the $I_{3}$ fiber and one fishtail fiber, call it
$F_{1}$, together with the section $E_{3}$. We blow up the double
point of $F_{1}$ and denote the exceptional sphere by $E_{10}$. We
further blow up at three points of the $I_{3}$ fiber, two of them on
the -2 sphere intersecting the section ($E_{12}$ and $E_{13}$) and
the third at the intersection of the remaining -2 spheres
($E_{11}$). Finally, we smooth the transverse intersection of
$F_{1}$ with $E_{3}$. The above procedure, the outcome of which is
depicted in Figure \ref{3}, provides the desired embedding.
\end{proof}

\begin{figure}[h] \centering
\input{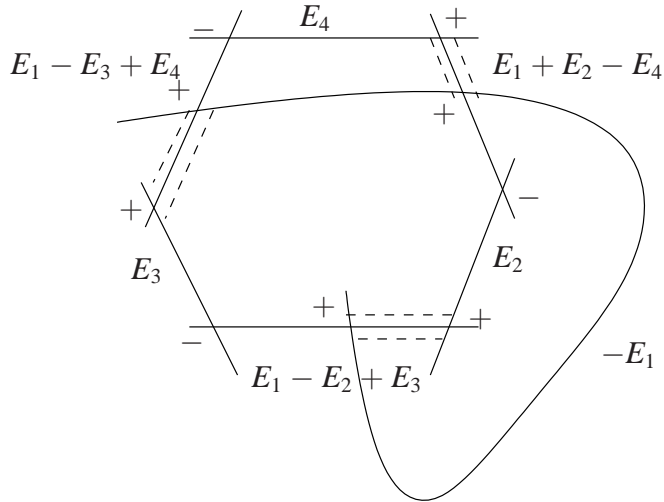}
\caption{$P_{1}$ in $\sharp 4\overline{\mathbb{CP} ^{2}}$} \label{4}
\end{figure}

\begin{remark}
\label{rem:2}
It is not hard to see that $Y_{1}$ bounds a rational
homology ball $B_{1}$. To this end, we can use the fact that $P_{1}$
embeds in $\sharp 4\overline{\mathbb{CP} ^{2}}$ (See Figure
\ref{4}). The closure of the complement of this embedding with
reversed orientation is a rational ball. Alternatively, we can use
\cite{Ne} and construct such a rational ball explicitly.
\end{remark}

\begin{theorem}
\label{thm:X1'}
$X_{1}'=(\mathbb{CP}^{2} \sharp 13
\overline{\mathbb{CP} ^{2}}-int(P_{1}))\bigcup_{Y_{1}}B_{1}$ is
homeomorphic to $\mathbb{CP}^{2} \sharp 9 \overline{\mathbb{CP}
^{2}}$.
\end{theorem}
\begin{proof} First, we will prove that $X_{1}'$ is simply
connected. The embedding of $P_{1}$ into $\mathbb{CP}^{2} \sharp 13
\overline{\mathbb{CP} ^{2}}$ constructed above has a simply
connected complement since the circles in the boundary $Y_{1}$ of
the complement are homotopically trivial in the complement(here we
are using the fact that rational surfaces are simply connected). In
fact, it suffices to prove this for the normal circles to the -3
framed spheres $C_{2}-E_{1}-E_{2}-\sum_{4}^{9}E_{i}-2E_{10}$ and
$3H-2E_{1}-E_{2}-\sum_{4}^{9}E_{i}-E_{11}$ and this can be done
easily using disks in one fishtail fiber and $E_{7}$ respectively.
In addition, the map $\pi_{1}(\partial B_{1}) \to \pi_{1}(B_{1})$
induced by the natural embedding is surjective, as one can see
applying Van Kampen's Theorem for the decomposition $\sharp 4
\overline{\mathbb{CP} ^{2}}=P_{1}\bigcup _{Y_{1}} ( \sharp 4
\overline{\mathbb{CP} ^{2}} \backslash P_{1})=P_{1} \bigcup
_{\partial B_{1}} \overline{B_{1}} $, where $\overline{B_{1}}$
denotes $B_{1}$ with opposite orientation. Thus, $X_{1}'$ is indeed
simply connected. Now the statement of the theorem follows from
Freedman's Theorem (\cite{Fr}), after computing the Euler
characteristic and the signature of the two manifolds.
\end{proof}

\subsection{An exotic smooth structure on $ \mathbb{CP}^{2} \sharp
8 \overline{\mathbb{CP} ^{2}}$.}
For our second construction,
we will make use of the existence result stated in Proposition 3.

\begin{prop}
\label{prop:3}
There is an elliptic fibration of E(1) $ \to
\mathbb{CP}^{1}$ with a singular $I_{5}$ fiber, 7 fishtail fibers
and one section.
\end{prop}

\begin{flushleft}
See the appendix for a brief discussion of this.
\end{flushleft}

Let $P_{2}$ denote the plumbing tree of spheres \\

\noindent \setlength{\unitlength}{3cm}
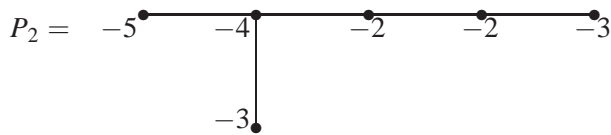
        \begin{picture}(3,.5)
          \put(.6,.5){\circle*{.05}}
          \put(.4,.4){$-5$}
          \put(1.1,.5){\circle*{.05}}
          \put(.9,.4){$-4$}
          \put(1.6,.5){\circle*{.05}}
          \put(1.5,.4){$-2$}
          \put(2.1,.5){\circle*{.05}}
          \put(2,.4){$-2$}
          \put(2.6,.5){\circle*{.05}}
          \put(2.5,.4){$-3$}
          \put(1.1,0){\circle*{.05}}
          \put(.9,0){$-3$}
          \put(.6,.5){\line(1,0){2}}
          \put(1.1,0){\line(0,1){.5}}
          \put(0,.4){$P_{2}=$}
         \end{picture}\\

\noindent and call $Y_{2}$ the boundary of $P_{2}$. It is not hard
to see that

\begin{prop}
\label{prop:4}
$P_{2}$ embeds into $\mathbb{CP}^{2} \sharp 14
\overline{\mathbb{CP} ^{2}}$.
\end{prop}

\begin{proof} Consider the fibration of Proposition \ref{prop:3}.
Blow up two double points in two of the fishtail fibers. Then
perform three further blow-ups at the $I_{5}$ fiber, one at the
intersection point of two -2 spheres and the other two on the -2
sphere intersecting the section. Finally, smooth out the
intersections of the two fishtails with the section to get the
desired embedding.
\end{proof}

\begin{remark}
\label{rmk:3}
To prove that $P_{2}$ bounds a rational homology ball,
we will once again use an appropriate embedding of this manifold,
i.e. the embedding of $P_{2}$ in $\sharp 6\overline{\mathbb{CP}
^{2}}$.
\end{remark}

\begin{theorem}
\label{thm:X2'}
$X_{2}'=(\mathbb{CP}^{2} \sharp 14
\overline{\mathbb{CP} ^{2}}-int(P_{2}))\bigcup_{Y_{2}}B_{2}$ is
homeomorphic to $\mathbb{CP}^{2} \sharp 8 \overline{\mathbb{CP}
^{2}}$.
\end{theorem}
\begin{proof} Simple connectivity of the complement of $P_{2}$ in
$\mathbb{CP}^{2} \sharp 14 \overline{\mathbb{CP} ^{2}}$ can be
proven using disks in a fishtail fiber and some of the exceptional
spheres. In addition to that, surjectivity of the map
$\pi_{1}(\partial B_{2}) \rightarrow \pi_{1}(B_{2})$ induced by the
natural embedding follows from an application of Van Kampen's
theorem for $\sharp 6\overline{\mathbb{CP} ^{2}}$, completely
analogous to the one in Claim 2 of Proof of Theorem \ref{thm:X1'}.
These facts, together with Freedman's classification theorem, lead
to the proof of the theorem.
\end{proof}

We would here like to point out that the plumbing trees of spheres
we have used so far, i.e. $P_{1}$ and $P_{2}$, belong in the
category of manifolds studied by Neumann in \cite{Ne}. In his
notation $P_{1} \sim M(0;(1,1),(3,2),(3,2),(3,2))$ (p=q=r=2) and
$P_{2} \sim M(0;(1,1),(3,2),(5,4),(5,2))$ (p=q=2,r=4). In addition,
note that using Neumann's results in this paper, one can find
rational balls bounded by $ \partial P_{i}, i = 1,2, $
with known handlebody decompositions. \\

For our next two constructions, we will need to combine the
techniques used above with knot surgery in a double node
neighborhood, as it was introduced in \cite{FiSt} by R. Fintushel
and R. Stern.

\subsection{An exotic smooth structure on
$ \mathbb{CP}^{2} \sharp 7 \overline{\mathbb{CP} ^{2}}$.}

\begin{figure} [h]
\centering
\input{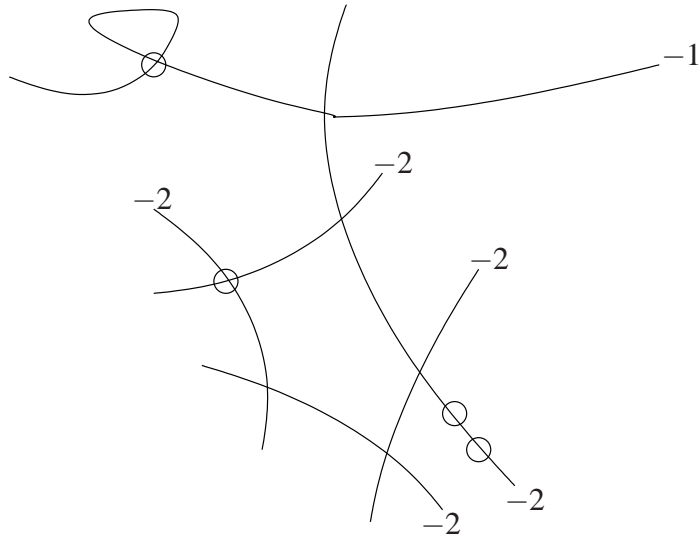}
\caption{The 4 blow-ups performed in construction 3} \label{5}
\end{figure}

Here, we will use the fibration for E(1) described in
construction 2. We will also use a double node neighborhood D
containing two of the fishtails of our fibration which have the same
monodromy and we will perform knot surgery along a regular fiber in
this neighborhood with a knot K having the properties listed in
\cite{FiSt}. The result of knot surgery will be to remove a smaller
disk from the disk $E_{5}$, which is a section in our original
picture, and to replace it with the Seifert surface of K. This will
give us a pseudo-section, that is a disk with a positive double
point in $H_{2}(D_{k}, \partial ; \mathbb{Z})$ and, after 4 blow-ups
as indicated in Figure \ref{5}, an embedding of $P_{2}$ in $
\mathbb{CP}^{2} \sharp 13 \overline{\mathbb{CP} ^{2}} $ (see Figure
\ref{6}). Our claim is that after blowing down we will get a
manifold $X_{3}'$ homeomorphic to $\mathbb{CP}^{2} \sharp 7
\overline{\mathbb{CP} ^{2}}$. Using disks on a fishtail fiber and on
the exceptional spheres $E_{10}$ and $E_{12}$, one can prove that
$\pi_{1}(\mathbb{CP}^{2} \sharp 13 \overline{\mathbb{CP} ^{2}} -
P_{2})$ is trivial. The rest of the argument is very similar to the
proofs of theorems 2 and 3 above and is therefore left as an
exercise for the reader. \\

\textbf{Another example ($X_{3}''$)}

The reader may have already noticed that it could be possible to
construct an exotic smooth structure on $ \mathbb{CP}^{2} \sharp 7
\overline{\mathbb{CP} ^{2}}$ by blowing down the 
next
plumbing tree in our family, that is $M(0;(1,1),(3,2),(7,6),(11,2))$
(p=q=2,r=6) in Neumann's notation, inside $\mathbb{CP}^{2} \sharp 15
\overline{\mathbb{CP} ^{2}}$.

This alternative construction can indeed be carried out using a
fibration of E(1) with an $I_{7}$ fiber and five fishtail fibers. We
give a short outline of this construction by indicating that the
double points of three of the fishtail fibers will be blown up so
that after smoothing out the intersection points of these fibers
with the section we can get a -7 sphere. One of the two remaining
fishtail fibers will then be used in proving simple connectivity of
the complement of the plumbing tree in $\mathbb{CP}^{2} \sharp 15
\overline{\mathbb{CP} ^{2}}$.

\begin{figure}  [h]
\centering
\input{4.pstex_t}
\caption{$P_{2}$ in $ \mathbb{CP}^{2} \sharp 13
\overline{\mathbb{CP} ^{2}} $} \label{6}
\end{figure}

\begin{remark}
\label{rmk:4}
A question most naturally arising here is whether
$X_{3}''$ is a member of the family of $X_{3}'$'s or not.
\end{remark}

\subsection{An exotic smooth structure on $ \mathbb{CP}^{2} \sharp
6 \overline{\mathbb{CP} ^{2}}$.}

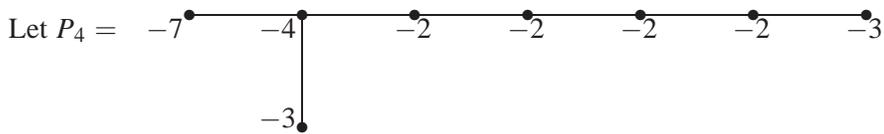
        \setlength{\unitlength}{3cm}
         \begin{picture}(4,.5)
          \put(.8,.5){\circle*{.05}}
          \put(.6,.4){$-7$}
          \put(1.3,.5){\circle*{.05}}
          \put(1.1,.4){$-4$}
          \put(1.8,.5){\circle*{.05}}
          \put(1.7,.4){$-2$}
          \put(2.3,.5){\circle*{.05}}
          \put(2.2,.4){$-2$}
          \put(2.8,.5){\circle*{.05}}
          \put(2.7,.4){$-2$}
          \put(3.3,.5){\circle*{.05}}
          \put(3.2,.4){$-2$}
          \put(3.8,.5){\circle*{.05}}
          \put(3.7,.4){$-3$}
          \put(1.3,0){\circle*{.05}}
          \put(1.1,0){$-3$}
          \put(.8,.5){\line(1,0){3}}
          \put(1.3,0){\line(0,1){.5}}
          \put(0,.4){Let $P_{4}=$}
         \end{picture}\\
\\

At this point, we only briefly note that starting with a fibration
of E(1) with one $I_{7}$ fiber and five fishtail fibers and
performing knot surgery along a regular fiber in a double node
neighborhood (as in 3.3) together with five appropriate blow-ups we
get an embedding of $P_{4}$ in $\mathbb{CP}^{2} \sharp 14
\overline{\mathbb{CP} ^{2}}$. Along the lines of our previous
arguments, it can be proven that $P_{4}$ bounds a rational homology
ball and that blowing down along $\partial P_{4}$ gives a manifold
$X_{4}'$ homeomorphic to $ \mathbb{CP}^{2} \sharp 6
\overline{\mathbb{CP} ^{2}}$.

\section{Computation of Seiberg--Witten invariants}

\begin{theorem}
\label{thm:4}
There is a characteristic cohomology class $K_{i}' \in
H^{2}(X_{i}')$ with $SW_{X_{i}'}(K_{i}') \neq 0.$
\end{theorem}

\begin{corollary}
\label{cor:1}
The 4--manifold $X_{i}'$ is not diffeomorphic to
$\mathbb{CP}^{2} \sharp (10-i) \overline{\mathbb{CP}^{2}}$.
\end{corollary}

\begin{proof} It is known that the SW invariants of $\mathbb{CP}^{2}
\sharp (10-i) \overline{\mathbb{CP}^{2}}, i \in \{1,2,3,4\}$, are
trivial, because of the existence of a metric with positive scalar
curvature. This, combined with Theorem \ref{thm:4} and the fact that
the SW invariants are diffeomorphism invariants, leads to a proof of
the corollary.
\end{proof}

\begin{proof}[Proof of \fullref{thm:4}]
We will apply Theorem \ref{thm:main} to all four constructions.
\begin{itemize}
\item[ (i) ] 1st construction: \\
Denote $\mathbb{CP}^{2} \sharp 13 \overline{\mathbb{CP} ^{2}}$ by
$X_{1}$. $Y_{1}$ is a monopole L--space (see \cite{KMOS} for a
proof) and $\vert H_{1} (Y_{1}) \vert = 81$. In addition, $B_{1}$
and $P_{1}$ are negative definite 4--manifolds. Consider
\begin{displaymath}
K_{1} \in H^{2}(X_{1};\mathbb{Z}), K_{1}(H)=3, K_{1}(E_{i})=1, i \in
{1,2,...,13}
\end{displaymath}
and denote $K_{1} \vert _{X_{1}- int(P_{1})}$ by $K_{1 \vert}$.
$K_{1 \vert}$ extends as a characteristic cohomology class to
$X_{1}'$ (proved using the embedding of $P_{1}$ in $\sharp 4
\overline{\mathbb{CP}^{2}}$ and more specifically that $K_{1}$
evaluates on the spheres of $P_{1}$ in the same way that the
canonical class of $\sharp 4 \overline{\mathbb{CP}^{2}}$ evaluates
on them). Denote this extension of $K_{1 \vert}$ by $K_{1}'$.
Finally, consider
\begin{equation}\label{eq:a1}
a_{1}=6H-2E_{1}-2E_{2}- \sum_{i=4}^{9}2E_{i}-E_{10}-E_{12}-E_{13}.
\end{equation}
Note that for such an $a_{1}$, the following conditions hold : \\
$a_{1} \in H_{2}(\mathbb{CP}^{2} \sharp 13 \overline{\mathbb{CP}
^{2}}; \mathbb{Z}), a_{1} \cdot a_{1} \ge 0, H \cdot a_{1} > 0,
K_{1}(a_{1})<0$, $a_{1}$ is represented in ($\mathbb{CP}^{2} \sharp
13 \overline{\mathbb{CP} ^{2}} -
int(Y_{1})$). \\
Now, since $H \cdot a_{1} > 0$, $K_{1}(a_{1}) < 0$ implies the
existence of a wall between PD(H) and $PD(a_{1})$. In the chamber
corresponding to PD(H), $SW_{X_{1}}$ is trivial, since we have a
positive scalar curvature metric. The wall crossing formula implies
therefore that
$SW_{X_{1},a_{1}}(K_{1})= \pm 1$. \\
By the dimension formula (2), $d(K_{1})=0$ and $d(K_{1}')=0$ since d
remains unchanged by our operation. It follows from the preceding
analysis that we can apply Theorem \ref{thm:main} to our case and
thus get that
\begin{displaymath}
SW_{X_{1}',a_{1}}(K_{1}')= \pm 1,
\end{displaymath}
which completes the proof of Theorem \ref{thm:X2'} for our first construction. \\

Note that there is no ambiguity about the chambers in the blown down
manifold, since the wall crossing formula combined with the
dimension formula for SW invariants implies that for a 4--manifold M
with $b_{2}^{+}(M)=1$ and $b_{2}^{-}(M) \leq 9 $ there is only one
chamber.

\item[ (ii) ] 2nd construction:
Asking for the analogous to the above conditions to be fullfilled,
one can easily compute that
\begin{equation}\label{eq:a2}
a_{2} = 7h - 3e_{1}-2 \sum_{2}^{9}e_{i}-e_{12}-e_{13}-2e_{14}
\end{equation}
is a cohomology class that can be used for the computation of the
Seiberg--Witten invariants in this case.

\item[ (iii)-(iv) ] In a similar fashion, one can carry out the
computations for the remaining cases.
\end{itemize}
\end{proof}

\begin{note}
D. Gay and A. Stipsicz recently proved in \cite{GaSt} that Wahl type
diagrams provide examples of plumbing trees that can be
symplectically blown down. Making use of their results, it follows
immediately that the manifolds $X_{1}', X_{2}'$ and $X_{3}''$
constructed above are symplectic.
\end{note}

\section{Appendix: elliptic fibrations of E(1)}
We give explicit constructions for some of the elliptic
fibrations $E(1) \to \mathbb{CP}^{1}$ used in the paper. Note that
the existence of such fibrations can also be verified using the
monodromies of the singular fibers.

\subsection{ A fibration of E(1) with a singular fiber of type
$I_{3}$ and nine fishtail fibers. }
\begin{figure} [ht]
\centering
\input{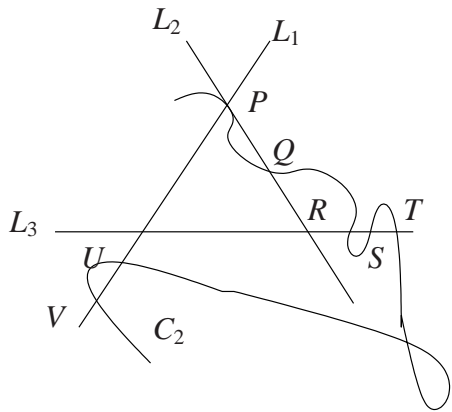}
\caption{Curves generating the pencil} \label{7}
\end{figure}

Let $C_{1}$ and $C_{2}$ be
two complex curves in the complex projective plane, such that :
$C_{2} = \{ [x:y:z] \in \mathbb{CP}^{2} \vert
p_{2}(x,y,z)=x^{3}+zx^{2}-zy^{2}=0 \} $ or any curve isotopic to
this so that it will give rise to a fishtail fiber in $E(1)$,
$C_{1}$ is the union of three lines - $L_{1}, L_{2},L_{3}$ - defined
by an equation of the form $p_{1}(x,y,z)=0$ and $C_{1},C_{2}$
intersect as indicated in Figure \ref{7}.

\begin{figure} [ht]
\centering
\input{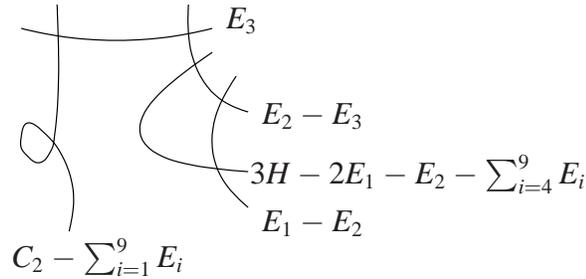}
\caption{A section, a fishtail and an $I_{3}$ fiber in E(1)}
\label{8}
\end{figure}

Perform three infinitely close blow-ups at the base point P and six
further blow-ups at the base points Q,R,S,T,U,V. After doing so, the
curves $C_{1}$ and $C_{2}$ get locally separated and the pencil of
elliptic curves
$C_{t}=C_{[t_{1}:t_{2}]}=\{(t_{1}p_{1}+t_{2}p_{2})^{-1}(0)\},
[t_{1}:t_{2}] \in \mathbb{CP}^{1}$ provides a well defined map $
\mathbb{CP}^{2} \sharp 9 \overline{\mathbb{CP} ^{2}} \to
\mathbb{CP}^{1}$. It is easy to check that starting with the curves
$C_{1},C_{2}$ and performing the nine blow-ups as described above,
one can detect an $I_{3}$ singular fiber, a fishtail fiber and
$E_{3}$ as a section (See also figure \ref{8}). Using the equations
defining $C_{1}$ and $C_{2}$, the remaining singular fibers can be
determined as well.

\subsection{A fibration of E(1) with a singular $I_{5}$ $(I_{7})$ fiber
and seven (five) fishtail fibers.}
Such fibrations can also be
easily constructed starting with appropriate generating curves for
the pencil.

Even more easily, one can verify the existence of such fibrations
combinatorially, using their monodromy. We will do so for one of the
two cases.

According to \cite{Ko}, the monodromies of our singular fibers are
as follows : \\
The monodromy of the fishtail fiber $I_{1}$ is $a = \left( \begin{array}{cc} 1 & 1 \\
0 & 1
  \end{array} \right)$
and the monodromy of the singular fiber $I_{k}$ is $ \left( \begin{array}{cc} 1 & k \\
0 & 1
\end{array}
  \right)$. If $ b = \left( \begin{array}{cc} 1 & 0 \\ -1  & 1
  \end{array} \right)$, then $b = (ab)a(ab)^{-1}$, and so b also represents the monodromy of a fishtail fiber.

In addition, $(a^{3}b)^{3} = I  \Leftrightarrow
a^{5}(a^{-2}ba^{2})abaaab= I $, which means that the fibration with
an $I_{5}$ and seven fishtail fibers over the disk extends to a
fibration over $S^{2}$. The classification of genus-1 Lefschetz
fibrations shows that the resulting fibration is an elliptic
fibration on E(1).

\end{document}